\documentclass[12pt]{amsart}
\usepackage{amsfonts,amssymb,amsmath,amsthm}

\newtheorem{thm}{Theorem}

\newtheorem{prop}[thm]{Proposition}

\begin{document}

\title[ZC for Small Groups]{The status of the Zassenhaus conjecture for small groups} 

\author[B\"achle]{Andreas B\"achle$^*$} 
\address{Vakgroep Wiskunde, Vrije Universiteit Brussel, Pleinlaan 2,
1050 Brussels, Belgium}
\email{abachle@vub.ac.be}

\author[Herman]{Allen Herman$^\dag$}
\address{Department of Mathematics and Statistics, University of Regina, Regina, Canada S4S 0A2} 
\email{Allen.Herman@uregina.ca}

\author[Konovalov]{Alexander Konovalov}
\address{School of Computer Science, University of St. Andrews, North Haugh, St. Andrews, Fife, KY16 SX9, U.K.}
\email{alexk@mcs.st-andrews.ac.uk}

\author[Margolis]{Leo Margolis$^\ddag$}
\address{Departamento de matem\'aticas, Facultad de matem\'aticas, Universidad de Murcia, 30100 Murcia, Spain}
\email{leo.margolis@um.es}

\author[Singh]{Gurmail Singh} 
\address{Department of Mathematics and Statistics, University of Regina, Regina, Canada S4S 0A2}
\email{singh28g@uregina.ca}

\thanks{$^*$ This author is a postdoctoral researcher of the FWO (Research Foundation Flanders).}
\thanks{$^\dag$ This author's research has been supported by an NSERC Discovery Grant.} 
\thanks{$^\ddag$ This author's research has been supported by a Marie Curie grant from EU project 705112-ZC} 
\date{August 31, 2016}



\begin{abstract}
We identify all small groups of order up to $288$ in the GAP Library for which the Zassenhaus conjecture on rational conjugacy of units of finite order in the integral group ring cannot be established by an existing method.  The groups must first survive all theoretical sieves and all known restrictions on partial augmentations (the HeLP$^+$ method).  Then two new computational methods for verifying the Zassenhaus conjecture are applied to the unresolved cases, which we call the quotient method and the partially central unit construction method.  To the cases that remain we attempt an assortment of special arguments available for units of certain orders and the lattice method.   In the end, the Zassenhaus conjecture is verified for all groups of order less than 144 and we give a list of all remaining cases among groups of orders 144 to 287. 
\end{abstract} 

\maketitle

\medskip
The Zassenhaus conjecture for torsion units of integral group rings states that any torsion unit of $\mathbb{Z}G$ should be conjugate in $\mathbb{C}G$ to a unit of $\pm G$.  We will abbreviate the conjecture as ZC1.  Though this conjecture was made in the 1960s, it has only been established for a few families of groups.  In this article, we study the status of the conjecture for groups of small order.  This effort follows an earlier one by H\"ofert and Kimmerle that dealt with groups of order up to 71 \cite{HK2006}.  Using recent theoretical and computational advances, we establish ZC1 for all groups of order less than 144, and give an accurate account of the groups of order 144 to 287 for which ZC1 cannot be established by the methods currently at our disposal (see Table \ref{Unresolved}.)  

\medskip
The first step is to apply theoretical sieves.  First we eliminate all nilpotent groups \cite{W91}, cyclic-by-abelian groups \cite{CMdR13}, and groups whose derived subgroup is contained in a Sylow subgroup for some prime $p$ \cite[Theorem 1.2]{Hertw06}.  Groups of the form $H \times C_2$ for which $\mathbb{Z}H$ satisfies ZC1 can also be eliminated \cite[Corollary 3.3]{HK2006}, since we are inducting on group order.   A further type of metabelian group described in \cite{MRSW87} can also be sieved.  

Let $\varepsilon$ denote the augmentation map on $\mathbb{C}G$.  A unit $u$ of $\mathbb{Z}G$ is normalized if $\varepsilon(u)=1$.  We will denote the normalized units of $\mathbb{Z}G$ by $\mathrm{V}(\mathbb{Z}G)$.  The partial augmentation of $u = \sum_g u_g g$ on the conjugacy class $x^G$ of $x \in G$ is $\varepsilon_x(u) = \sum_{y \in x^G} u_y$.  An important connection between torsion units in $\mathrm{V}(\mathbb{Z}G)$, established in \cite[Theorem 2.5]{MRSW87}, says that the torsion units satisfying the property:

\smallskip
(*) for all $n \ge 0$, there is only one conjugacy class $x^G$ in $G$ such that $\varepsilon_x(u^n) \ne 0$ 

\smallskip
\noindent are precisely the torsion units that are conjugate in $\mathbb{C}G$ to normalized trivial units.  This result has inspired an approach that investigates the potential existence of normalized torsion units with more than one nonzero partial augmentation.  The Luthar-Passi method \cite{LP}, its modular extension by Hertweck \cite{Hertw-PA}, the Cohn-Livingstone and power congruence criteria (see \cite[Remark 6]{Bovdi-Hertweck}), and Hertweck's criteria that $\varepsilon_x(u) \ne 0$ implies the order of $x$ must divide the order of $u$ \cite[Theorem 2.3]{Hertw-PA}  produce several restrictions on partial augmentations that a normalized torsion unit of $\mathbb{Z}G$ must satisfy.  The first and fourth authors recently produced the GAP package {\tt help} \cite{help} that implements this collection of methods, to which we refer as HeLP$^+$.   

After the sieve, the HeLP$^+$-method is applied to the remaining list of groups.  For this step we use the GAP package \cite{help}, which utilizes the integer optimization programs {\tt normaliz} \cite{normaliz} and {\tt 4ti2} \cite{4ti2}.  Each HeLP$^+$-solution consists of the lists of partial augmentations for a possible normalized torsion unit and its powers.  A nontrivial HeLP$^+$-solution is one for which $\varepsilon_x(u) \ne 0$ for more than one conjugacy class $x^G$ of $G$.  For each nontrivial HeLP$^+$-solution, we have to show that there is no torsion unit of $\mathrm{V}(\mathbb{Z}G)$ that produces these partial augmentations.  As this is an inductive approach, we can also eliminate a nontrivial HeLP$^+$ solution for $G$ by mapping it modulo a proper normal subgroup $N$ to a nontrivial list of partial augmentations summing in $G/N$.  If a normalized unit $u$ with these partial augmentations exists in $\mathbb{Z}G$, then its image would be a normalized unit with nontrivial partial augmentations in $\mathbb{Z}[G/N]$, whose existence would already be ruled out.   Elimination by this ``quotient method'' is quite useful in groups that have several normal subgroups.  

The HeLP package also gives the spectrum for the torsion unit under each irreducible representation of the group.   This information is equivalent to the partial augmentation information.  It is used directly in the partially central construction and in the lattice method.  In the partially central construction, we directly check if the central components of $u$ in the Wedderburn decomposition of $\mathbb{C}G$ lie in the image of $\mathbb{Z}G$.  If they do not, the unit cannot be conjugate to an element of $\mathbb{Z}G$.  The partially central construction and the quotient method are implemented with software (see \cite{github}).  

To the groups with nontrivial HeLP$^+$-solutions that remain after these computational procedures, we check a few more additional criteria available for units of particular orders by hand, and finally check if the lattice method of \cite{BM15} can be applied.  In the end, only a handful of groups of order between 144 and 287 survive all of these tests (see Table \ref{Unresolved}.) 

Our notation is based on GAP's.  We will refer to the group identified by {\tt SmallGroup(n,r)} in GAP's Small Groups Library as $SG(n,r)$.  For brevity we only record one nontrivial HeLP$^+$-solution in each $\operatorname{Aut}(G)$-orbit.  We denote conjugacy classes with their labels in GAP, and also use this notation for partial augmentations, so for example, ($4a,4c,4f)=(2,1,-2)$ means the partial augmentations of the unit would be $2$ on the class of elements of order $4$ GAP labels with $4a$ by GAP, $1$ on the class $4c$, $-2$ on the class $4f$, and $0$ on all other conjugacy classes. 

\begin{table}[h]\centering
\caption{Unresolved nontrivial HeLP$^+$-solutions among groups of order up to 287.}\label{Unresolved}
\begin{tabular}{c|c|l} 
\mbox{GAP Id} & \mbox{Structure} & \mbox{unresolved cases: PAs[powers]}  \\ \hline
(144,117) & $C_3^2 \rtimes D_{16}$ & (6b,6d,6e)=(-1,1,1)[2b,3b] \\
(144,119) & $C_3^2 \rtimes Q_{16}$ & (12a,12b,12d)=(1,1,-1)[2a,3a,4a,6a]  \\
(150,5) & $C_5^2 \rtimes S_3$ & (10a,10c)=(2,-1)[2a,5a] \\
(160,234) & $C_2^4 \rtimes D_{10}$ & (2a,2b,2c)=(1,-1,1) \\
(192,955) & $C_2^4 \rtimes D_{12}$ & (2a,2c,2f)=(1,1,-1),(1,-1,1) \\
(192,973) & $C_2^3 \rtimes SG(24,8)$ & (8a,8bc,8d) = (1,-1,1)[2b,4db] \\
(192,974) & $C_2^3 \rtimes  SG(24,8)$ &(8a,8bc,8d) = (1,-1,1)[2c,4a] \\
(192,975) & $C_2^3 : SG(24,8)$ & (8a,8bc,8d) = (1,-1,1)[2a,4fc] \\
(192,976) & $C_2^3 : SG(24,8)$ & (8a,8bc,8d) = (1,-1,1)[2b,4eb] \\
(192,1489) & $C_2^3 : S_4$ & (8a,8b,8c)=(-1,1,1)[2b,4e] \\
(192,1490) & $C_2^3 : S_4$ & (8a,8b,8c) = (-1,1,1)[2c,4d] \\
 (216,33) & $C_3^3 \rtimes Q_8$ & (12a,12c,12f) = (1,1,-1)[2a,3a,4a,6a] \\
& & \qquad \qquad \qquad (0,2,-1)[2a,3c,4c,6c] \\
 (216,35) & $C_3^3 \rtimes D_8$ & (12a,12b)=(2,-1)[2c,3c,4d,6e] \\ 
 (216,37) & $C_3^3 \rtimes D_8$ & (6a,6e,6h) = (1,1,-1)[2a,3a], \\
 & & \qquad \qquad \qquad (0,2,-1)[2c,3c] \\
 (216,153) & $C_3^2 \rtimes SL(2,3)$ & (3a,3d,6a)= (1,-1,1),(-1,1,1)[3c,2a] \\
 (240,91) & $A_5 \rtimes C_4$ & (4a,4c,4d)=(1,0,-1,1)[2a] \\
\end{tabular}
\end{table}

\noindent {\bf Remarks}.  (i). In Table \ref{Unresolved}, the nontrivial lists of partial augmentations of the unit $u$ and the conjugacy classes in $G$ of powers $u^d$ for $d$ dividing the order $o(u)$ of $u$ are given.  Since we only seek a counterexample of minimal unit order, we assume $u^d$ is always conjugate to an element of $G$ for any proper divisor $d$ of $o(u)$.    We use an order convention for multiple nontrivial solutions; so $(8a,8bc,8d) = (1,-1,1)[2a,4fc]$ in Table \ref{Unresolved} dentoes the information for the two solutions $(8a,8b,8d)=(1,-1,1)[2b,4f]$ and  $(8a,8c,8d) = (1,-1,1)[2a,4c]$.

\noindent (ii). $SG(24,8) \simeq C_3 \rtimes D_8$ is one of three groups with this structure description. 

\noindent (iii). $SG(240,91)$ is a (central) quotient of $GL(2,5)$.  ZC1 was verified for $GL(2,5)$ by Bovdi and Hertweck in \cite{Bovdi-Hertweck}.  However, it does not follow from their arguments that ZC1 holds for $G$.  

\bigskip
\noindent {\bf 1.  Groups of order less than 288 that survive our sieves.}  

\bigskip
We abbreviate the group identified as {\tt SmallGroup(n,r)} in GAP's Small Groups Library as $SG(n,r)$.  Our first theoretical sieve removes all nilpotent, cyclic-by-abelian, and $p$-by-abelian groups, and groups satisfying case (a) of \cite[Main Theorem]{MRSW87}.   Since we are looking for a minimal counterexample, we can also eliminate groups of the form $C_2 \times H$, where the group $H$ is known to satisfy ZC1.   After this sieve, $1121$ groups of order less than 288 remain. 

To these groups we apply the HeLP$^+$ method.  This leaves a list of 99 groups.    

Next, we apply the quotient method.  If the list of partial augmentations in a nontrivial HeLP$^+$-solution maps onto a nontrivial list of partial augmentations in any proper quotient of the group, it is eliminated by our inductive argument.  This leaves a list of 35 groups, which appear in Table \ref{AfterQM}. 

\begin{table}[h] \centering
\caption{Cases remaining after sieves, HeLP$^+$, and the quotient method. }\label{AfterQM}
\begin{tabular}{l|l}
\mbox{ Order } & \mbox{ GAP Id \#, (unit order)}  \\ \hline 
48 & 30(4) \\
72 & 40(6) \\
96 & 65(8), 186(4), 227(2) \\
144 & 117(6), 119(12), 182(6) \\
150 & 5(10) \\
160 & 234(2,4) \\
168 & 43(6) \\
192 & 182(4), 186(16), 955(2), 958(4,8) \\
        & 973(8), 974(8), 975(8), 976(8), 987(4), 990(4), \\
				& 1473(4), 1478(4), 1489(8), 1490(8),  \\
				& 1491(4), 1492(4), 1493(2) \\
200 & 43(10) \\
216 & 33(12), 35(12), 37(6), 153 (3,6), 161 (3) \\  
240 & 91(4) \\ 
\end{tabular}
\end{table} 

\bigskip
\noindent {\bf Remark}.  The smallest groups our methods do not cover turn out to have order 144, so it is at order 288 where the quotient method becomes ineffective.   

\bigskip
\noindent {\bf 2.  Partially central unit constructions} 

\bigskip
The second, third, and fifth authors have developed a computer program that can construct torsion units of $\mathbb{C}G$ that have the partial augmentations and powers in a given HeLP-solution (see \cite{github} and \cite{HS15}).  Each nontrivial HeLP$^+$-solution has accompanying spectral information.  A complete set of explicit irreducible representations of the group is used to construct a unit of $\mathbb{C}G$ with the desired spectrum, and it follows it is conjugate in $\mathbb{C}G$ to any torsion unit with the same partial augmentations and powers.  We then take $e$ to be the sum of all the centrally primitive idempotents of $\mathbb{C}G$ for which $ue$ is central in $\mathbb{C}Ge$.  The centrality can be observed directly from the spectral information since $u$ must represent diagonally under these irreducible representations.  Finally, we determine a subset of $Ge$ that gives a $\mathbb{Z}$-basis of $\mathbb{Z}Ge$, and express $ue$ in this basis.  If $ue \not\in \mathbb{Z}Ge$, then no $\mathbb{C}G$-conjugate of $u$ will ever be an element of $\mathbb{Z}G$, and so this HeLP$^+$-solution will not produce a counterexample to ZC1.  When this occurs, we simply say that the unit $u$ is {\it partially central}.  

When constructing units, we write $e_{3abc}$ for the sum of the centrally primitive idempotents of $\mathbb{C}G$ corresponding to $\chi_{3a}$, $\chi_{3b}$, and $\chi_{3c}$.  If we write $e_{124}$, then this means the sum of all centrally primitive idempotents of $\mathbb{C}G$ corresponding to all irreducible characters of degrees $1$, $2$, and $4$.  We use, for example, $g_{4b}$ to denote a group element in the conjugacy class of $G$ labelled $4b$ by GAP, $\chi_{3d}$ for the $4$th irreducible character of degree $3$ in GAP. 

\medskip
{\bf SG(48,30), order $4$.}  There is one $\operatorname{Aut}(G)$-orbit of nontrivial HeLP$^+$-solutions that survive the quotient method.  Its partial augmentations are  $(4a,4c,4d) = (-1,1,1)$, and its square lies in the class $2a$. Its spectral information shows the central part of the unit is $ue_{13cd}$.  

{\small
$$\begin{array}{r|ccccccc}
& \chi_{1abcd} & \chi_{2a} & \chi_{2b} & \chi_{3a} & \chi_{3b} & \chi_{3c} & \chi_{3d} \\ \hline
spec(\mathfrak{X}(u)) & 1,-1,i,-i & (i,-i) & (1,-1) & (1,1,1) & (-1, -1, -1 ) & (i,i,-i) & (i,-i,-i) \\
\end{array} $$
}

Using this spectral information and explicit irreducible representations constructed in GAP for each $\chi \in \operatorname{Irr}(G)$, we can construct a representative of $u$ in $\mathbb{C}G$.   With our notational conventions, the $u$ our program constructs looks like 
$$ u = g_{4b} e_{123cd} +g_{2a}e_{3a}-g_{2a}e_{3b}. $$ 
With our program we then construct an explicit $\mathbb{Z}$-basis of $\mathbb{Z}Ge_{13cd}$ and write $u e_{13cd}$ in terms of this basis.  When we do this the coefficients we get are non-integer rational numbers, so this means $u$ is partially central.   So the partially central unit construction completes the verification of ZC1 for this group.    

\medskip
As we have implemented this partially central check with our software, we just give a summary of the cases where it works.  

\medskip
{\bf SG(96,65), order $8$:}   There is one $\operatorname{Aut}(G)$-orbit of nontrivial HeLP$^+$-solutions that survive the quotient method:  $(8a,8d,8g) = (1,-1,1)$, with $u^2 \in 4a$, $u^4 \in 2a$.   The central part of the unit is $ue_{13efgh}$ and it does not lie in $\mathbb{Z}Ge_{13efgh}$.  So ZC1 holds for $G$. 

\medskip
{\bf SG(96,186), order $4$:}  There is one $\operatorname{Aut}(G)$-orbit of solutions left to resolve:  $(4b,4g,4h) = (-1,1,1)$, with $u^2 \in 2b$.  Our calculations show its central part $ue_{13efgh} \not\in \mathbb{Z}Ge_{13efgh}$.  So ZC1 holds for $G$. 

\medskip
{\bf SG(96,227), order $2$:}  One orbit of solutions to resolve: $(2b,2c,2d)=(1,1,-1)$.  The unit we construct is partially central.  So ZC1 holds for $G$. 

\medskip
{\bf SG(192,182), order $4$:}  Only one solution is left to resolve: $(4c,4d,4g)=(1,1,-1)$, $u^2 \in 2a$.  Since $ue_{13gh} \not\in \mathbb{Z}Ge_{13gh}$, this verifies ZC1 for $G$. 

\medskip
{\bf SG(192,186), order $16$:}  One solution is left to resolve: $(16a,16e,16k)=(1,-1,1)$, $u^2 \in 8a$, $u^4 \in 4a$, and $u^8 \in 2a$.   Since $ue_{13ijk\ell mnop} \not\in \mathbb{Z}Ge_{13ijk\ell mnop}$, ZC1 holds for $G$. 

\medskip
{\bf SG(192,955), order $2$:}  There are five $\operatorname{Aut}(G)$-orbits of order $2$ solutions left to resolve:
$$\begin{array}{rcl}
(2a,2b,2c,2d,2f)&=&(1,0,-1,0,1), (0,1,-1,0,1), \\
&  & (0,0,1,1,-1),(1,0,1,0,-1),(0,0,2,0,-1).
\end{array}$$
For the first, second, and fourth of these, the central part of the unit we construct for them lies in $\mathbb{Z}Ge_{}$, so the partially central method does not resolve them.  For the third and fifth, the central part of the unit does not lie in $\mathbb{Z}Ge$, so these cases are resolved.  The other three nontrivial HeLP$^+$-solutions for this group remain unresolved. 

\medskip
{\bf SG(192,958), order $4$:}  One solution left to resolve: $(4b,4g,4j)=(1,-1,1)$, $u^2 \in 2b$.  The unit we construct for this solution is partially central, so this verifies ZC1 for units of order $4$. 

\medskip
{\bf SG(192,958), order $8$:}  One orbit of solutions left to resolve: $(8b,8h,8n)=(1,-1,1)$, with $u^2 \in 4a$, and $u^4 \in 2b$.  $ue_{13abcdghij} \not\in \mathbb{Z}Ge_{13abcdghij}$, so the unit is partially central.  This verifies ZC1 for $G$. 

\medskip
{\bf SG(192,987), order $4$:} There are four $\operatorname{Aut}(G)$-orbits of HeLP$^+$-solutions remaining to be resolved, which we can summarize as $(4abcd,4e,4f)=(1,1,-1)$, $u^2 \in 2a$. Our construction shows all four of these units are partially central.  So this verifies ZC1 for $G$.  

\medskip
{\bf SG(192,990), order $4$:} There are three $\operatorname{Aut}(G)$-orbits of nontrivial solutions to consider: $(4abc,4d,4e)=(1,1,-1)$, with $u^2 \in 2c$.  All three are shown to be partially central by our program.  So ZC1 is verified for $G$. 

\medskip
{\bf SG(192,1473), order $4$:} One orbit of solutions to resolve: $(4b,4h,4k)=(1,-1,1)$, $u^2 \in 2c$, which is shown to be partially central by our construction.  So this resolves ZC1 for $G$. 

\medskip
{\bf SG(192,1478), order $4$:} One orbit of solutions to resolve: $(4c,4i,4m)=(1,-1,1)$, $u^2 \in 2b$, which is shown to be partially central by our construction.  So this resolves ZC1 for $G$. 

\medskip
{\bf SG(192,1491), order $4$:} One orbit of solutions to resolve: $(4a,4d,4e)=(1,-1,1)$, $u^2 \in 2c$.  The unit our program constructs for this solution is partially central, so this verifies ZC1 for $G$. 

\medskip
{\bf SG(192,1492), order $4$:} There are two orbits of solutions remaining: $(4ab,4c,4e)=(1,-1,1)$, with $u^2 \in 2b$.  The units our program constructs for these solutions are partially central, so this verifies ZC1 for $G$. 

\medskip
{\bf SG(192,1493), order $2$:} There is one orbit of solutions to resolve, represented by $(2a,2b,2e)=(-1,1,1)$.  The unit our program gives for this solution is partially central, so this verifies ZC1 for $G$. 

\bigskip
\noindent {\bf 3.  Methods for units of special order.}  

\bigskip
The first of the special results we can apply is Proposition 4.2 of \cite{Hertw06} concerning $p$-adic conjugacy of units.  An easy integral consequence of this $p$-adic result is the following one, relevant for the $p$-subgroup version of the Zassenhaus conjecture ($p$-ZC3):    

\begin{prop} \label{pZC3}
Let $N$ be a normal $p$-subgroup of $G$.  Suppose $U$ is a finite subgroup of $\mathrm{V}(\mathbb{Z}G)$ which maps to $1$ under the natural map modulo $N$.  Then $U$ is conjugate in $\mathbb{Q}G$ to a subgroup of $N$. 
\end{prop} 

{\bf SG(160,234), order $2$:} There are two orbits of solutions to consider: $(2a,2b,2c)=(1,1,-1)$ and $(2b,2c,2d)=(1,1,-1)$.  For the latter, the classes $2b$, $2c$, and $2d$ are contained in a normal subgroup $N$ of order $16$, so we can apply Proposition \ref{pZC3} to resolve it.  The other solution remains unresolved. 

\medskip
{\bf SG(192,955), order $2$:}  The three nontrivial solutions remaining to be resolved are $(2a,2b,2c,2f) = (1,0,-1,1),(1,0,1,-1),(0,1,-1,1)$.  The group has a normal subgroup $N$ of order $32$ that contains the classes $2b$, $2c$, and $2f$.  An application of Proposition \ref{pZC3} resolves the third of these solutions.  The first two remain unresolved. 

The next lemma we can use is \cite[Proposition 2]{Hertw08b}:  

\begin{prop} \label{vanishingppart}
Let $N$ be a normal $p$-subgroup of $G$.  Suppose $u$ is a torsion unit of $\mathbb{Z}G$ with augmentation $1$ whose image modulo $N$ has strictly smaller order than $u$.  Then $\varepsilon_g(u)=0$ for every $g \in G$ whose $p$-part has strictly smaller order than the $p$-part of $u$. 
\end{prop}

{\bf SG(72,40), order $6$:}  The nontrivial HeLP$^+$ solutions lie in two distinct $\operatorname{Aut}(G)$-orbits, whose partial augmentations are: $(2a,2c,6a)=(1,1,-1),(-1,1,1)$.  The group has a normal subgroup $N$ of order $9$, and the image of $u$ modulo $N$ in both cases would have order $2$.  Both of these cases fail the criteria of Proposition \ref{vanishingppart} at the prime $p = 3$. 

\medskip
{\bf SG(144,182), order $6$:}  The same reasoning applies.  There are two orbits of nontrivial solutions, with $(2a,2b,6a)=(1,1,-1),(-1,1,1)$.  There is a normal subgroup $N$ of order $9$ for which $u$ mod $N$ has order $2$.   So these partial augmentations of $u$ also fail Proposition \ref{vanishingppart}.

\medskip
{\bf SG(160,234), order $4$:}  There is one orbit of nontrivial solutions to resolve: 
$(2a,4a)=(2,-1)$, $u^2 \in 2c$.   As noted earlier, $G$ has a normal subgroup of order $16$ containing $g_{2c}$ but not $g_{2a}$.  Modulo $N$, $u$ will have order $2$,  so Proposition \ref{vanishingppart} applies.  So normalized torsion units of order $4$ in $\mathbb{Z}G$ will be rationally conjugate to elements of $G$.  

\medskip
{\bf SG(168,43), order $6$:}  There are ten nontrivial HeLP$^+$-solutions of order $6$: 
$$\begin{array}{rcl} 
(3a,3b,6a,6b) &=& (1,2,-1,-1),(1,-1,-1,2),(2,1,-1,-1),(-1,1,2,-1),(-1,1,1,0), \\
                         &=& (1,-1,0,1),(2,1,-2,0),(1,2,0,-2),(3,0,-2,0),(0,3,0,-2).
												\end{array}$$
The group $G$ has a nontrivial normal subgroup $N$ of order $8$ generated by the class $2a$.  Since $|G/N|=21$, any torsion unit $u$ of $\mathbb{Z}G$ that has order $6$ will be mapped modulo $N$ to an element of order $3$.  By applying Proposition \ref{vanishingppart} with $p=2$ we can eliminate every nontrivial solution with $\varepsilon_{3a}(u) \ne 0$ or $\varepsilon_{3b}(u)\ne 0$.  So this verifies ZC1 for $G$. 

\medskip
{\bf SG(200,43), order $10$:}  There are two $\operatorname{Aut}(G)$-orbits of nontrivial HeLP$^+$-solutions to resolve: $(2a,2c,10a)=(-1,1,1),(1,1,-1)$.  $G$ has a normal subgroup of order $25$, so Proposition \ref{vanishingppart} tells us the partial augmentations of a normalized torsion unit of $\mathbb{Z}G$ with order $10$ should vanish on elements of order $2$.  So these solutions are resolved, and this verifies ZC1 for $G$.  

We will need another result \cite[Theorem 5.3]{DJ96} where ZC3 for $p$-subgroups is known. 

\begin{prop}\label{lcc-pZC3}
Let $G$ be a finite solvable group, and suppose $L$ is the last nontrivial term of the lower central series of $G$.  If $p$ is a prime dividing $|L|$ for which $p^4$ does not divide $|G|$, then any finite $p$-subgroup of normalized torsion units is rationally conjugate to a subgroup of $G$. 
\end{prop}

Of course, the last nontrivial term of the lower central series is the smallest normal subgroup $L$ of $G$ for which $G/L$ is nilpotent. 

\medskip
{\bf SG(216,153), order $3$:}  For this group, if $L$ is the last non-trivial term of the lower central series of $G$, then $|L|=72$.  By Proposition \ref{lcc-pZC3}, normalized units of $\mathbb{Z}G$ with $2$- or $3$-power order are rationally conjugate to elements of $G$.  

\medskip
{\bf SG(216,161), order $3$:}  For this group, if $L$ is the last non-trivial term of the lower central series of $G$, then $|L|=27$.  By Proposition \ref{lcc-pZC3}, normalized units of $\mathbb{Z}G$ with $3$-power order are rationally conjugate to elements of $G$.  This verifies ZC1 for this group. 

\bigskip
\noindent {\bf 5. Applying the Lattice Method.} 

\bigskip
Let $G = SG(216, 153)$, the special affine group $\mathbb{F}_3^2 \rtimes \operatorname{SL}(2,3)$. After the arguments above there remain four $\operatorname{Aut}(G)$-orbits of possibly non-trivial partial augmentations for units of order $6$ in $\mathrm{V}(\mathbb{Z}G)$. Two of these we will exclude using the lattice method introduced in \cite{BM15}. We will use parts of the character table of $G$ and the decomposition matrix of $G$ given in table \ref{CT72}. The values of all characters on classes not given in the table are integral. 

\begin{table}[h]\centering
\caption{Character table of $G$}\label{CT72}
\begin{tabular}{l c c c c c c}\hline
{} & 1a & 2a & 3a & 3c & 3d & 6a \\ \hline\hline
$\chi_{1a}$ & 1 & 1 & 1 & 1 & 1 & 1 \\
$\chi_{2a}$ & 2 & -2 & -1 & -1 & -1 & 1 \\
$\chi_{3a}$ & 3 & 3 & 0 & 0 & 0 & 0 \\
$\chi_{8a}$ & 8 & 0 & 2 & 2 & -1 & 0 \\
\end{tabular}
\end{table} 

{\begin{table}[h]\centering
\caption{Decomposition matrix of $G$ for the prime $3$}\label{DM72}
\begin{tabular}{l c c c}\hline
{} & $\varphi_{1a}$ & $\varphi_{2a}$ & $\varphi_{3a}$ \\ \hline\hline
$\chi_{1a}$ & 1 & 0 & 0  \\
$\chi_{2a}$ & 0 & 1 & 0  \\
$\chi_{3a}$ & 0 & 0 & 1  \\
$\chi_{8a}$ & 1 & 2 & 1 \\
\end{tabular}
\end{table} 
 
$2a$ is the only class of involutions in $G$ and in both possibilities of units we are going to study $u^2$ is rationally conjugate to elements in $3c$. Denote by $\zeta$ a primitive complex 3rd root of unity. Assume first that $(\varepsilon_{3a}(u), \varepsilon_{3d}(u), \varepsilon_{6a}(u)) = (-2,2,1).$ Then for representations $D_{2a}$ and $D_{8a}$ corresponding to the characters $\chi_{2a}$ and $\chi_{8a}$ respectively we find the following eigenvalues for $u$.
\[D_{2a}(u) \sim {\rm{diag}}(-\zeta, -\zeta^2), \quad D_{8a}(u)  \sim {\rm{diag}}(\zeta, \zeta^2, \zeta,\zeta^2,-1,-1,-1,-1). \] 
By a theorem of Fong \cite[Corollary 10.13]{Isaacs} we may assume that the representations $D_{2a}$ and $D_{8a}$ are realized over a $3$-adically complete discrete valuation ring which is unramified over the $3$-adic integers. Denote by $\bar{.}$ the reduction modulo the maximal ideal of $R$, also with respect to modules, and let $k$ be the residue class field of $R$. Let $L_{2a}$ and $L_{8a}$ be $RG$-lattices corresponding to $D_{2a}$ and $D_{8a}$ respectively. Then when viewed as $k\langle \bar{u} \rangle$-modules we have by \cite[Proposition 1.3]{BM15} that $\bar{L}_{2a} \cong \bar{L}_{2a}^+ \oplus \bar{L}_{2a}^-$ and $\bar{L}_{8a} \cong \bar{L}_{8a}^+ \oplus \bar{L}_{8a}^-$ are such that all the composition factors of $\bar{L}_{2a}^+$ and $\bar{L}_{8a}^+$ are trivial while the composition factors of $\bar{L}_{2a}^-$ and $\bar{L}_{8a}^-$ are non-trivial. Moreover by \cite[Propositions 1.3, 1.4]{BM15} we know that $\bar{L}_{2a}$ is a $2$-dimensional indecomposable module while $\bar{L}_{8a}^-$ is the direct sum of four $1$-dimensional summands.  But since $\bar{L}_{2a}^-$ is a sub- or factor module of $\bar{L}_{8a}^-$ by the decomposition numbers in table \ref{DM72}, this contradicts the existence of $u$.

Next assume $(\varepsilon_{3a}(u), \varepsilon_{3d}(u), \varepsilon_{6a}(u)) = (2,-2,1).$ We will use similar notation as in the paragraph above. Then
\[D_{3a}(u) \sim {\rm{diag}}(1,\zeta, \zeta^2), \quad D_{8a}(u)  \sim {\rm{diag}}(-\zeta, -\zeta^2, -\zeta, -\zeta^2,1,1,1,1). \] 
Then viewed as $k\langle \bar{u} \rangle$-modules we conclude, again by \cite[Propositions 1.3, 1.4]{BM15}, that $\bar{L}_{3a}^+$ contains a direct indecomposable summand of dimension at least $2$ while $\bar{L}_{8a}^+$ is the direct sum of four $1$-dimensional modules. This again contradicts the decomposition numbers given in table \ref{DM72}.

\end{document}